\documentstyle[12pt]{article}
\newcommand{\be}{\begin{equation}}
\newcommand{\ee}{\end{equation}}
\newcommand{\bea}{\begin{eqnarray}}
\newcommand{\ena}{\end{eqnarray}}
\newcommand{\bae}{\begin{eqnarray}}
\newcommand{\eae}{\end{eqnarray}}
\newcommand{\sect}[1]{\setcounter{equation}{0}\section{#1}}

\newcommand{\cC}{{\cal C}}

\newcommand{\cG}{{\cal G}}
\newcommand{\cH}{{\cal H}}
\newcommand{\cI}{{\cal I}}

\newcommand{\cO}{{\cal O}}

\newcommand{\cU}{{\cal U}}
\newcommand{\cZ}{{\cal Z}}
\newcommand{\CC}{\rlap {\raise 0.4ex \hbox{$\scriptstyle |$}}
      \hskip -0.2em \hbox{{\rm C}}}
\newcommand{\ZZ}{\hbox {\sf Z\hskip-5pt Z}}
\newcommand{\diag}{\, \hbox{diag} \,}
\newcommand{\mod}{\, \hbox{mod} \,}
\newcommand{\id}{\, \hbox{id} \,}
\newcommand{\Tr}{\, \hbox{Tr}}

\newcommand{\uqs}{\mbox{${\cal U}_{qs}$}}
\newcommand{\Aqs}{{\cal A}_{qs}}
\def\uLpm{\uL^{\pm}}
\def\uLp{\uL^{+}}
\def\uLm{\uL^{-}}
\newcommand{\uA}{\underline{A}}
\newcommand{\uL}{\underline{L}}
\newcommand{\uS}{\underline{S}}
\newcommand{\us}{\underline{S}}
\newcommand{\uY}{\underline{Y}}
\newcommand{\uDelta}{\underline{\Delta}}
\newcommand{\ueps}{\underline{\epsilon}}
\newcommand{\utimes}{\underline{\otimes}}

\newcommand{\wh}[1]{\widehat{#1}}

\newcommand{\ip}[2]{\left\langle \Big. {#1},{#2} \right\rangle}
\newcommand{\bip}[2]{\left\langle\!\!\left\langle \Big. {#1},{#2}
                     \right\rangle\!\!\right\rangle}
\newcommand{\tip}[2]{\left \lfloor \! \! \left \lfloor \Big. {#1},{#2}
                     \right \rfloor \! \! \right \rfloor}
\font\petit=cmr7
\font\moyen=cmr10
\font\grand=cmr12
\font\Grand=cmr12 scaled \magstep1
\newcommand{\enslapp}{{\petit E}{\moyen N}{\grand S}{\Grand L}{\grand
    A}{\moyen P}{\petit P}}

\def\hRot{\hat{R}_{12}}

\def\hRtth{\hat{R}_{23}}

\def\qit{q^{-2}}
\def\Lpm{L^{\pm}}
\def\Lp{L^{+}}
\def\Lm{L^{-}}
\def\ie{\hbox{\it i.e.}}

\def\fe{& = &}

\def\ff{\nn \\}

\def\uot{\underline{\otimes}}
\def\hR{\hat{R}}
\def\men#1{~(\ref{#1})}
\def\mens#1{~\ref{#1}}
\def\nn{\nonumber}
\def\uU{1}
\def\ra{\rightarrow}

\def\ot{\otimes}

\begin{document}
\renewcommand{\thefootnote}{\fnsymbol{footnote}}
\newpage
\pagestyle{empty}
\setcounter{page}{0}
\newcommand{\norm}[1]{{\protect\normalsize{#1}}}
\newcommand{\LAP}
{{\small E}\norm{N}{\large S}{\Large L}{\large A}\norm{P}{\small P}}
\newcommand{\sLAP}{{\scriptsize E}{\footnotesize{N}}{\small S}{\norm L}$
${\small A}{\footnotesize{P}}{\scriptsize P}}
\begin{minipage}{5.2cm}
\begin{center}
{\bf Groupe d'Annecy\\
\ \\
Laboratoire d'Annecy-le-Vieux de Physique des Particules}
\end{center}
\end{minipage}
\hfill
\hfill
\begin{minipage}{4.2cm}
\begin{center}
{\bf Groupe de Lyon\\
\ \\
Ecole Normale Sup\'erieure de Lyon}
\end{center}
\end{minipage}
\centerline{\rule{12cm}{.42mm}}

%

\vfill
\vfill
\begin{center}

  {\LARGE {\bf {\sf Classical and Quantum $sl(1|2)$ Superalgebras,
        Casimir Operators and Quantum Chain Hamiltonians}}}\\[1cm]

  \vfill

  {\large Daniel Arnaudon$^{a}$,
    Chryssomalis Chryssomalakos$^{b}$
    and Luc Frappat$^{c}$}

  \vfill

  {\em
    Laboratoire de Physique Th\'eorique }\LAP\footnote{URA 14-36
    du CNRS, associ\'ee \`a l'Ecole Normale Sup\'erieure de Lyon et \`a
    l'Universit\'e de Savoie.

    \indent
    $^a$ arnaudon@lapp.in2p3.fr

    \indent
    $^b$ chryss@lapp.in2p3.fr

    \indent
    $^c$ frappat@lapp.in2p3.fr


    }\\
    {\em Chemin de Bellevue BP 110, F-74941 Annecy-le-Vieux Cedex,
      France.}

\end{center}

\vfill

\begin{abstract}
  We examine in this paper the two parameter deformed superalgebra
  $\cU_{qs}(sl(1|2))$ and use the results in the construction of
  quantum chain Hamiltonians. This study is done both in the framework of
  the Serre presentation and in the $R$-matrix scheme
  of Faddeev, Reshetikhin and Takhtajan (FRT).
  We show that there exists an
  infinite number of Casimir operators, indexed by integers $p \ge 2$
  in the undeformed case and by $p \in \ZZ$ in the deformed
  case, which obey quadratic relations. The construction of the dual
  superalgebra of functions on $SL_{qs}(1|2)$ is also given and
  higher tensor product representations are discussed.
  Finally, we construct quantum chain Hamiltonians based on the Casimir
  operators.
  In the deformed case we find two Hamiltonians which describe
  deformed $t-J$ models.
\end{abstract}

\vfill
\vfill

\rightline{\enslapp-A-505/95}
\rightline{q-alg/9503021}
\rightline{March 95}

\newpage
\pagestyle{plain}

\sect{Introduction \label{Intro}}

\indent

A Lie algebra $\cG$ being given, it is known that the coproduct of the
elements of the centre $C(\cG)$ of $\cG$, evaluated in a
representation of $\cG$, gives two-site quantum chain Hamiltonians,
which are invariant under the Lie algebra $\cG$ under consideration.
A physically important case is the $t-J$ model at the supersymmetric
point which can be obtained by the above procedure using the second order
Casimir operator and the fundamental representation
of the simple Lie superalgebra $sl(1|2)$. Recently, the
authors of ref. \cite{FeiYue} found a $q$-deformation of the $t-J$ model,
invariant under the quantum group $U_q(sl(1|2))$.

However, there are reasons to wonder whether the Hamiltonians one
can construct in this way are unique.

A simple Lie superalgebra being given, there exist many inequivalent
simple root bases (\ie\  many Cartan matrices) \cite{KacSuper} which
can be used
to describe the superalgebra in the Serre presentation
\cite{ScheuSerre}.
Unlike in the classical case, where the Hopf structure is essentially
unique, the quantum deformation of $\cU(sl(1|2))$ exhibits the novel
feature that the Hopf structures that are natural in each basis are
not related by the transformation law between the corresponding bases.
One can thus suspect that they lead to unequal Hamiltonians.

On the other hand, different Casimir operators may
also lead, through the above procedure, to different Hamiltonians.

It is our purpose, in this paper, to explore these
possibilities. Accomplishing this necessitates an efficient formalism
for dealing with quantum superalgebras. There are already in the
literature approaches to an FRT-type~\cite{FRT} formulation of
quantum superalgebras~\cite{CPR,MRP} which, however, deal only with
universal enveloping superalgebras (see also~\cite{LS,PSVdJ,SSVZ}
for alternative approaches from a different perspective). As our
study calls for representations of elements of the above algebra
(or its tensor powers), we give a new FRT construction of $\uqs
(sl(1|2))$, complete with its dual superalgebra of functions on
$SL_{qs}(1|2)$.

\indent

The paper is organised as follows.
In Sect. \ref{Class} we recall the structure  of the classical
superalgebra $sl(1|2)$,
using the Serre presentation in the fermionic
basis. We show that the centre of its universal enveloping
superalgebra $\cU(sl(1|2))$ contains an infinite set of Casimir
operators, indexed by integers $p\ge 2$, that satisfy quadratic relations.
The quantum analogue $\cU_{qs}(sl(1|2))$ of $\cU(sl(1|2))$, with two
deformation parameters, is considered in Sect. \ref{Quantum I}. The
corresponding deformed Casimir operators, indexed by $p\in \ZZ$ are
given.  They also obey quadratic relations.
In Sect. \ref{quantum II} we recast $\uqs (sl(1|2))$ in its FRT
form, introducing a suitable $R$-matrix and matrices of generators
$\uLpm$. We also discuss available choices in the construction of
the dual Hopf algebra. A bosonised basis for the algebra is given
in Sect. \ref{Boson} and the representation of the coproduct of an
infinite set of Casimirs is computed.
Finally, in Sect. \ref{Hamil}, the above results are used to
construct three-state (per site) quantum chain Hamiltonians. The
$sl(1|2)$-invariant Hamiltonian is unique (up to the identity and a
normalization factor) and describes the well-known $t-J$ model at the
supersymmetric point.
A similar construction in the case of
$\uqs(sl(1|2))$-invariant Hamiltonians provides two
$qs$-deformations of the supersymmetric $t-J$ model.
In the case of open boundary conditions, an equivalence
transformation eliminates all but the usual deformation parameter $q$,
leading to two Hamiltonians, one being
the $(1,2)$ Perk--Schultz Hamiltonian \cite{PerkS}.
The results are finally summarised in Sect. \ref{Conclusion}.

\sect{The classical superalgebra $\cU(sl(1|2))$ \label{Class}}

\indent

Before studying the quantum case, we would like to recall some
features and do some comments about the (undeformed) simple Lie
superalgebra $sl(1|2)$. A classical simple Lie superalgebra can be
described in many inequivalent bases (\ie\ there
exist many inequivalent simple root systems), each
one being associated to a particular Dynkin diagram.
The different possible Dynkin diagrams can be found by applying
generalised Weyl transformations associated to the fermionic roots of
the diagrams, until no new diagram appears. In
the case of $sl(1|2)$, the two possibilities are the
so-called distinguished and fermionic bases.

\subsection{Presentation in the fermionic basis}

\indent

Denoting by $H_i \ (i=1,2)$ and $E^{\pm}_i \ (i=1,2,3)$ the generators
of $sl(1|2)$, the $\ZZ_2$-gradation in the fermionic basis is such that
$H_1$, $H_2$ and $E^{\pm}_3$ are even elements and
$E^{\pm}_1$, $E^{\pm}_2$ are odd elements, that is
\be
\deg H_1 = \deg H_2 = \deg E^{\pm}_3 = 0 \ \ \ \mbox{and} \ \ \
\deg E^{\pm}_1 = \deg E^{\pm}_2 = 1
\label{Grading}
\ee
where $\deg X$ stands for the degree of the generator $X$.

In the Serre presentation, the commutation relations
in the fermionic basis read
\bea
&&[H_i, H_j] = 0 \ \ \ \mbox{and} \ \ \
[H_i, E^{\pm}_j] = \pm a_{ij} E^{\pm}_j \;, \nonumber \\
&&\{E^{\pm}_1, E^{\pm}_1\} = \{E^{\pm}_2, E^{\pm}_2\} =
\{E^{\pm}_1, E^{\mp}_2\} = 0 \;, \nonumber \\
&& \{E^+_1, E^-_1\} = H_1 \;, \nonumber \\
&& \{E^+_2, E^-_2\} = H_2 \;;
\label{CommClass}
\ena
$(a_{ij})$ is the Cartan matrix in the fermionic basis:
\be
(a_{ij}) = \left(\begin{array}{cc} 0 & -1 \\ -1 & 0 \end{array}
\right) \;.
\label{CartanMat}
\ee
The generators $E^{\pm}_3$ are defined by
\be
E^{\pm}_3 = \{ E^{\pm}_1 , E^{\pm}_2 \}
\label{E3}
\ee
and the Serre relations can be written as
\be
[E^{\pm}_1 , E^{\pm}_3] = [E^{\pm}_2 , E^{\pm}_3] = 0 \;.
\label{SerreRelClass}
\ee
It follows from the previous equations that the remaining
relations among the generators read
\be
\begin{array}{ll}
\bigg. [E^+_3, E^-_1] = E^+_2 \;, \qquad
&[E^+_3, E^-_2] = E^+_1 \;, \\
\bigg. [E^-_3, E^+_1] = - E^-_2 \;, \qquad
&[E^-_3, E^+_2] = - E^-_1 \;, \\
\bigg. [E^+_3, E^-_3] = H_1 + H_2  \;. &
\end{array}
\label{AlgClass}
\ee

The universal enveloping superalgebra
$\cU \equiv \cU(sl(1|2))$ can be endowed with
a classical super-Hopf structure, the $\ZZ_2$-graded coproduct
$\uDelta$ being given by
\be
\uDelta(X) = X \utimes 1 + 1 \utimes X
\label{DeltaClass}
\ee
for $X \in sl(1|2)$ and extending super-multiplicatively to
the entire $\cU$.
Here and in the following, we use an underlined notation for the
graded structures. For example, $\utimes$ denotes the $\ZZ_2$-graded
tensor product satisfying
\be
(A \utimes B)(C \utimes D) = (-1)^{\deg B.\deg C} (AC \utimes BD) \;.
\label{GradDelta}
\ee

\subsection{Casimir operators of $\cU$}

\indent

In this subsection, we focuse our attention on the centre of $\cU$.
We can construct a (countable) infinite set of
Casimir operators $\cC^{cl}_p$ where $p$ is an integer $\ge 2$. The
explicit expression of $\cC^{cl}_p$ is
\bea
\cC^{cl}_p
    &=&  +  \; H_1 H_2 \Big( H_1-H_2 \Big)^{p-2}\nonumber \\
    &&   - \; E^-_1 E^+_1 \Big(H_2 (H_1-H_2)^{p-2}
                        + (1-H_2) (H_1-H_2+1)^{p-2} \Big)\nonumber \\
    &&   - \; E^-_2 E^+_2 \Big(H_1 (H_1-H_2)^{p-2}
                        + (1-H_1) (H_1-H_2-1)^{p-2} \Big)\nonumber \\
    &&   - \; E^-_3 E^+_3 \Big(H_1-H_2 \Big)^{p-2}\nonumber \\
    &&   + \; E^-_3 E^+_2 E^+_1
              \Big( (H_1-H_2)^{p-2}-(H_1-H_2+1)^{p-2} \Big)\nonumber \\
    &&   + \; E^-_2 E^-_1 E^+_3
              \Big( (H_1-H_2)^{p-2}-(H_1-H_2-1)^{p-2} \Big)\nonumber \\
    &&   + \; E^-_2 E^-_1 E^+_2 E^+_1 \Big((H_1-H_2+1)^{p-2}
                        + (H_1-H_2-1)^{p-2} - 2 (H_1-H_2)^{p-2}
                        \Big) \;. \nonumber \\
    &&
\label{CasClass}
\ena
They satisfy the following relations
\be
\cC^{cl}_{p_1} \cC^{cl}_{p_2} = \cC^{cl}_{p_3} \cC^{cl}_{p_4}
\ \ \ \mbox{if} \ \ \ p_1 + p_2 = p_3 + p_4
\quad \mbox{and} \quad p_i \ge 2 \;.
\label{CasClassRel}
\ee
Nevertheless, it is important to note that none of these Casimir operators can
be
written as a polynomial function of the others.
However, because of these relations, the eigenvalues of only two Casimir
operators are enough to characterise a finite-dimensional highest
weight irreducible representation of $sl(1|2)$.

\indent

Let $h$ be the projection
\be
h: \cU = \cU^0 \oplus \left( \cU \cU^+ + \cU^- \cU \right)
\longrightarrow   \cU^0
\ee
within the direct sum, where $\cU^+$ and $\cU^-$ are the
subalgebras of $\cU$ generated respectively by $E^+_i$
and $E^-_i$ ($i=1,2$), while $\cU^0$ is the unital subalgebra
generated by $1$ and $H_i$ ($i=1,2$).
\\
The restriction $\bar h$ of $h$ to the centre $\cZ_{\cU}$
of $\cU$  is an
algebra morphism onto the algebra ${\cU^0}^W$
of polynomials in the Cartan
generators $H_1$ and $H_2$
invariant under the action of the Weyl group, \ie\
($H_1 \leftrightarrow -H_2$)\footnote{Note that the translated Weyl
group
coincides here with the Weyl group. Indeed, denoting by
$\rho_0$ and $\rho_1$ the half sums of positive even and odd roots in
the fermionic basis,
the translation vector $\rho = \rho_0 - \rho_1$ vanishes for $sl(1|2)$.}.
$\bar h$ is the Harish--Chandra homomorphism
\cite{KacSuperChar}.  It is injective \cite{ScheuCas}.
{From} Eq. (\ref{CasClass}),
its image $\bar h \left( \cZ_{\cU} \right)$
is $\{1\}\cup \cI$, where $\cI$ is the ideal in ${\cU^0}^W$
generated by the product
$H_1 H_2$. This result is compatible with the fact that the fields of
fractions of ${\cU^0}^W$ and $\bar h \left( \cZ_{\cU} \right)$
coincide \cite{KacSuperChar}.

\indent

In reference \cite{JarvisCas}, the formula
\be
\cC_p^{(m,n)} = \hbox{STr} \left({\Big. \wh{E}\; }^p \right)
\label{CasJarvis}
\ee
is given for the Casimir operators of the special linear superalgebras
$sl(m|n)$. The matrix $\wh{E}$ is defined by
$\wh{E}^A_B=(-1)^{\deg B}E^A_B$ ($A,B=1,\cdots,m+n$) where
the matrix $E$ collects the set of generators of
$sl(m|n)$, with $\deg B=1$ for $B=1,\cdots,m$
and $\deg B =-1$ for $B=m+1,\cdots,m+n$. The authors of
\cite{JarvisCas} proved that, in any finite-dimensional highest
weight irreducible representation, the Casimir operators
$\cC_p^{(m,n)}$ with $p\ge m+n$ can be expressed in terms of
the previous ones $\cC_r^{(m,n)}$ with $r < p$.
\\
Our Casimir operator $\cC^{cl}_{p}$ of Eq. (\ref{CasClass}) is a
linear combination of $\cC_p^{(1,2)}$ and of a
polynomial function in the  $\cC_r^{(1,2)}$ with $r<p$. It is
actually {\it the} linear combination that eliminates all the higher
degree terms in the raising and lowering generators $E^{\pm}_i$ and
the image of which by the Harish--Chandra isomorphism is, up to a
power of $H_1-H_2$, a monomial of degree one in $H_1 H_2$.
It is worthwhile to notice that in the case of $sl(n)$, this
combination is precisely $0$ for $p > n$, expressing the fact that
$\cC_p^{n}$ for $p > n$ is equal to a polynomial function in the $\cC_r^{n}$,
$r\le n$.
In the case of $sl(1|2)$ however, as operators,
none of the Casimirs $\cC^{cl}_{p}$
can be written in terms of the others.

\subsection{Comments on the distinguished basis}

\indent

Let us end this section with some comments on the distinguished
basis. The generalised Weyl transformation
\be
\begin{array}{lll}
\bigg. h_1 = - H_1 - H_2    & h_2 = H_2 & \\
\bigg. e_1^\pm = \pm E_3^\pm \qquad \quad
     & e_2^\pm = E_2^\mp \qquad \quad
     & e_3^\pm = \pm E_1^\pm
\end{array}
\label{ferm-to-dist}
\ee
defines the generators of the distinguished basis in terms of those
of the fermionic basis.
The relations in the Serre--Chevalley
distinguished basis read
\bea
&&[h_i, h_j] = 0 \ \ \ \mbox{and} \ \ \
[h_i, e^{\pm}_j] = \pm a'_{ij} e^{\pm}_j \;, \nonumber \\
&&\{e^{\pm}_2, e^{\pm}_2\} = [e^{\pm}_1, e^{\mp}_2] = 0 \;, \nonumber \\
&& [e^+_1, e^-_1] = h_1 \;, \nonumber \\
&& \{e^+_2, e^-_2\} = h_2 \;,
\label{CommClassd}
\ena
where
$(a'_{ij})$ is the Cartan matrix in the distinguished basis:
\be
(a'_{ij}) = \left(\begin{array}{cc} 2 & -1 \\ -1 & 0 \end{array}
\right) \;.
\label{CartanMatd}
\ee
All the results of this section obtained in the fermionic basis can be
rewritten {\it mutatis mutandis} in the distinguished basis using the
transformation Eq. (\ref{ferm-to-dist}).

\sect{The quantum superalgebra $\cU_{qs}(sl(1|2))$ \label{Quantum I}}
\subsection{Presentation in the fermionic basis}

\indent

We consider now the two-parameter quantum universal
enveloping superalgebra
$\uqs \equiv \cU_{qs}(sl(1|2))$. It is the two-parametric
deformation of $\cU(sl(1|2))$, defined by the generators $H_i$,
$E^{\pm}_i$
$(i=1,2)$ and unit 1, and by the relations in the
Serre--Chevalley fermionic basis:
\bea
&&[H_i, H_j] = 0 \ \ \ \mbox{and} \ \ \
[H_i, E^{\pm}_j] = \pm a_{ij} E^{\pm}_j \;, \nonumber \\
&&\{E^{\pm}_1, E^{\pm}_1\} = \{E^{\pm}_2, E^{\pm}_2\} =
\{E^{\pm}_1, E^{\mp}_2\} = 0 \;, \nonumber \\
&& \{E^+_1, E^-_1\} = [H_1]_q \; s^{-H_1} \;, \nonumber \\
&& \{E^+_2, E^-_2\} = [H_2]_q \; s^{H_2} \;,
\label{CommQuant}
\ena
with $\displaystyle [x]_q \equiv \frac{q^x-q^{-x}}{q-q^{-1}}\;$
(note that for $s=1$, this is the standard definition of $\cU_{q}(sl(1|2))$).
\\
Defining
\bea
&& E^+_3 = q s^{-1} E^+_1 E^+_2 + E^+_2 E^+_1 \;,
\nonumber \\
&& E^-_3 = E^-_1 E^-_2 + q^{-1} s^{-1} E^-_2 E^-_1 \;,
\label{E3Q}
\ena
the quantum Serre relations read
\bea
&&E^{\pm}_1 E^{\pm}_3 - q^{-1} s^{\pm 1} E^{\pm}_3 E^{\pm}_1 = 0 \;, \nonumber
\\
&&E^{\pm}_2 E^{\pm}_3 - q\ s^{\mp 1} E^{\pm}_3 E^{\pm}_2 = 0 \;.
\label{SerreRelQuant}
\ena
$\uqs$ is equipped with the $\ZZ_2$-gradation given in
Eq. (\ref{Grading}) and is
endowed with the structure of a super Hopf algebra by defining the
coproduct $\uDelta : \uqs \longrightarrow \uqs \otimes \uqs$, the antipode
$\uS : \uqs \longrightarrow \uqs$ and the counit $\ueps : \uqs
\longrightarrow \CC$ as follows:
\bea
&&
\!\!
\begin{array}{l}
\uDelta(H_i) = H_i \utimes 1 + 1 \utimes H_i \;,\\
\uDelta(E^+_1) = E^+_1 \utimes 1 + q^{H_1} s^{-H_1} \utimes E^+_1 \;,\\
\uDelta(E^+_2) = E^+_2 \utimes 1 + q^{H_2} s^{H_2} \utimes E^+_2 \;,\\
\uDelta(E^-_1) = E^-_1 \utimes q^{-H_1} s^{-H_1} + 1 \utimes E^-_1 \;,\\
\uDelta(E^-_2) = E^-_2 \utimes q^{-H_2} s^{H_2} + 1 \utimes E^-_2
\;,
\label{Deltaferm}
\end{array}
\\
&& \nonumber \\
&&
\!\!
\begin{array}{ll}
\uS(H_i) = -H_i \;, & \\
\uS(E^+_1) = - q^{-H_1} s^{H_1} E^+_1\;, & \uS(E^+_2) = - q^{-H_2}
s^{-H_2} E^+_2 \;, \\
\uS(E^-_1) = -E^-_1 q^{H_1} s^{H_1}\;, & \uS(E^-_2) = -E^-_2 q^{H_2}
s^{-H_2}\;,
\end{array}
\\
&& \nonumber \\
&& \ueps(H_i) = \ueps(E^{\pm}_i) = 0 \;.
\label{Antipode}
\ena
Notice that the antipode in Eq. (\ref{Antipode}) is
super-antimultiplicative,
\ie
\be
\uS(AB) = (-1)^{\deg A. \deg B} \uS(B) \uS(A) \;.
\label{GradS}
\ee
We check by inspection that $S^2=\id$.
\\
It follows from Eqs. (\ref{CommQuant}) that the remaining relations among the
generators read
\bea
&&[E^+_3, E^-_1] = q s E^+_2 q^{-H_1} s^{-H_1} \;, \nonumber \\
&&[E^+_3, E^-_2] = E^+_1 q^{H_2} s^{H_2} \;, \nonumber \\
&&[E^-_3, E^+_1] = - E^-_2 q^{H_1} s^{-H_1} \;, \nonumber \\
&&[E^-_3, E^+_2] = -q^{-1} s E^-_1 q^{-H_2} s^{H_2} \;, \nonumber \\
&&[E^+_3, E^-_3] = [H_1+H_2]_q \; s^{H_2-H_1}
\label{AlgQuant}
\ena
with the associated Hopf structure
\bea
&& \uDelta(E^+_3) = E^+_3 \utimes 1 + \lambda s^{H_2} E^+_1 q^{H_2} \utimes
E^+_2
+ q^{H_2+H_1} s^{H_2-H_1} \utimes E^+_3 \;, \nonumber \\
&& \uDelta(E^-_3) = E^-_3 \utimes q^{-H_2-H_1} s^{H_2-H_1}
- \lambda E^-_2 \utimes q^{-H_2} E^-_1 s^{H_2} + 1 \utimes E^-_3 \;, \\
&& \nonumber \\
&& \uS(E^+_3) = -q^{-H_1-H_2} s^{H_1-H_2} E^+_3
+ \lambda q^{-H_1-H_2} s^{H_1-H_2-1} E^+_1 E^+_2 \;, \nonumber \\
&& \uS(E^-_3) = -E^-_3 q^{H_1+H_2} s^{H_1-H_2}
- \lambda E^-_2 E^-_1 q^{H_1+H_2} s^{H_1-H_2-1} \;, \\
&& \nonumber \\
&& \ueps(E^+_3) = \ueps(E^-_3) = 0 \;,
\label{E3Delta}
\ena
where $\lambda = q-q^{-1}$.

\subsection{Casimir operators of $\uqs$}

\indent

One can compute the Casimir operators of
$\uqs$. As in the undeformed case, one finds
by direct calculation a (countable) infinite set of Casimir operators
$\cC_p$, with $p$ here in $\ZZ$, given by
\bea
\cC_p = q^{(2p-1)(H_2-H_1)} & \bigg\{
& [H_1]_q [H_2]_q
+ E^-_1 E^+_1 s^{H_1-1} \left( q^{-2p+1} [H_2-1]_q - [H_2]_q
\right) \nonumber \\
&& + E^-_2 E^+_2 s^{-H_2-1} \left( q^{2p-1} [H_1-1]_q - [H_1]_q
\right)
- q^{-1} E^-_3 E^+_3 s^{H_1-H_2-1} \nonumber \\
&& \bigg. - q^{p-2} (q-q^{-1}) [p]_q E^-_2 E^-_1 E^+_3 s^{H_1-H_2-2}
\nonumber \\
&& \bigg. + q^{-p} (q-q^{-1}) [p-1]_q E^-_3 E^+_2 E^+_1 s^{H_1-H_2-1}
\nonumber \\
&& + q^{-1} (q-q^{-1})^2 [p]_q [p-1]_q E^-_2 E^-_1 E^+_2 E^+_1
s^{H_1-H_2-2} \bigg\} \;.
\label{CasQuant}
\ena
The Casimir operators $\cC_p$ satisfy the following relations
\be
\cC_{p_1} \cC_{p_2} = \cC_{p_3} \cC_{p_4} \ \ \ \mbox{if} \ \ \ p_1 +
p_2 = p_3 + p_4\;.
\label{CasQuantRel}
\ee
It is interesting to look at the classical limit of the $\cC_p$'s.
It is easy
to see that $\lim\limits_{q,s \rightarrow 1} \cC_p =
\cC^{cl}_2$ for all $p \in
\ZZ$. Moreover, the classical Casimir operators $\cC^{cl}_p$ with $p
\ge 3$ can be obtained as
limits of suitable linear combinations of the quantum Casimir
operators
$\cC_r$ involving coefficients with negative powers in $q-q^{-1}$ as
follows
\begin{eqnarray*}
\hbox{\hskip 4.1cm}
\cC^{cl}_p = \lim_{q,s \rightarrow 1} {1\over(q-q^{-1})^{p-2}}
\sum_{l=0}^{p-2} (-1)^l \left( {p-2 \atop l} \right) \cC_l \;.
\hbox{\hskip 3.4cm} (\simeq \pi)
\end{eqnarray*}

\indent

Let now $h$ be the projection
\addtocounter{equation}{1}
\be
h: \uqs = \uqs^0 \oplus \left( \uqs \uqs^+ + \uqs^- \uqs \right)
\longrightarrow   \uqs^0
\ee
within the direct sum, where $\uqs^0$, $\uqs^+$, $\uqs^-$ are the
subalgebras of $\uqs$ generated respectively by
$q^{\pm H_i}$ and $s^{\pm H_i}$, $E^+_i$, $E^-_i$ ($i=1,2$).
The restriction $\bar h$ of $h$ to the centre $\cZ_{\uqs}$
of $\uqs$  is an
algebra morphism onto the algebra of polynomials in the Cartan
generators $q^{\pm 2 H_1}$ and $q^{\pm 2 H_2}$,
invariant under the action of the Weyl group, \ie\
($H_1 \leftrightarrow - H_2$).
$\bar h$
is the quantum analogue of the
Harish--Chandra homomorphism.
{From} Eq. (\ref{CasQuant}),
its image $\bar h \left( \cZ_{\uqs} \right)$
is $\{1\}\cup \cI_{qs}$, where $\cI_{qs}$ is
the ideal generated by the product
$q^{H_1-H_2} [H_1]_q [H_2]_q$.

\subsection{Presentation in the distinguished basis}

\indent

As in the classical case, we can switch to the distinguished basis,
the change of basis being given by
\be
\begin{array}{lll}
\bigg. h_1 = - H_1 - H_2  \quad   & h_2 = H_2 & \\
\bigg. e_1^+ = s^{-1} E_3^+ &
       e_2^+ = E_2^- q^{-H_2} s^{-H_2} \quad &
       e_3^+ = s^{-1} E_1^+ \\
\bigg. e_1^- = - E_3^- &
       e_2^- = s^{-1} E_2^+ q^{H_2} s^{-H_2} \quad &
       e_3^- = - E_1^- \;.
\end{array}
\label{ferm-to-dist-qs}
\ee
This allows the derivation of the relations satisfied by
the generators in the distinguished basis.
The Casimir operators in the distinguished basis are hence provided from
(\ref{CasQuant}) using (\ref{ferm-to-dist-qs}).
\\
Unlike in the classical case, the natural Hopf structure in the
distinguished basis is not the one provided by this redefinition.
The natural coproduct of $\uqs$ in the distinguished basis is
actually given by (take $s=1$ to see that it is natural)
\bea
&&\tilde{\uDelta}(h_i) = h_i \utimes 1 + 1 \utimes h_i \;,\nonumber \\
&&\tilde{\uDelta}(e^+_1) = e^+_1 \utimes 1
                   + q^{h_1} s^{h_1+2h_2} \utimes e^+_1 \;,\nonumber \\
&&\tilde{\uDelta}(e^+_2) = e^+_2 \utimes 1
                   + q^{h_2} s^{-h_2} \utimes e^+_2 \;,\nonumber \\
&&\tilde{\uDelta}(e^-_1) = e^-_1 \utimes q^{-h_1} s^{h_1+2h_2}
                   + 1 \utimes e^-_1 \;,\nonumber \\
&&\tilde{\uDelta}(e^-_2) = e^-_2 \utimes q^{-h_2} s^{-h_2}
                   + 1 \utimes e^-_2
\;,
\label{Deltadist}
\ena
where now
\be
\deg h_1 = \deg h_2 = \deg e^{\pm}_1 = 0 \ \ \ \mbox{and} \ \ \
\deg e^{\pm}_2 = \deg e^{\pm}_3 = 1 \;.
\label{Gradingd}
\ee
The above can be considered as an alternative Hopf structure
compatible with the relations (\ref{CommQuant}) defining the fermionic
basis, in terms of which it reads
\bea
&&\tilde{\uDelta}(H_i) = H_i \utimes 1 + 1 \utimes H_i \;,\nonumber \\
&&\tilde{\uDelta}(E^+_1) = E^+_1 \utimes 1
                   + q^{-H_1} s^{-H_1} \utimes E^+_1
                   - \lambda q^{H_2} E^+_3 s^{-H_2} \utimes
                        q^{-H_2} E^-_2 s^{-H_2} \;,\nonumber \\
&&\tilde{\uDelta}(E^+_2) = E^+_2 \utimes q^{-2H_2}
                   + q^{-H_2} s^{H_2} \utimes E^+_2 \;,\nonumber \\
&&\tilde{\uDelta}(E^-_1) = E^-_1 \utimes q^{H_1} s^{-H_1}
                   + 1 \utimes E^-_1
                   - \lambda s^{-H_2} E^+_2 q^{H_2} \utimes
                        s^{-H_2} E^-_3  q^{-H_2} \;,\nonumber \\
&&\tilde{\uDelta}(E^-_2) = E^-_2 \utimes q^{H_2} s^{H_2}
                   + q^{2H_2} \utimes E^-_2 \;,\nonumber \\
&&\tilde{\uDelta}(E^+_3) = E^+_3 \utimes 1
                   + q^{-H_1-H_2} s^{H_2-H_1} \utimes E^+_3 \;,\nonumber \\
&&\tilde{\uDelta}(E^-_3) = E^-_3 \utimes q^{H_1+H_2} s^{H_2-H_1}
                   + 1 \utimes E^-_3 \;.
\label{Deltafermtilde}
\ena

In the construction of the Hamiltonians or for any other application, we
can equivalently use the distinguished basis with its natural Hopf
structure, or the fermionic basis with the new Hopf structure derived
from the latter.


\sect{Super-FRT Construction} \label{quantum II}

\indent

Anticipating the applications of section\mens{Hamil}, we attempt
now a formulation along the lines of the standard FRT construction.
Our motivation is twofold: on the one hand, we aim at a significant
computational simplification, due to the compactness of the
notation used in~\cite{FRT}, while, on the other hand, it is in this
formulation that our results can most easily be generalised to the
case of other quantum supergroups. We present our approach, as
already outlined in the introduction, in two stages. In this
section, continuing working in a $\ZZ_2$-graded environment, we find
the analogues of the $L^\pm$ matrices of~\cite{FRT} in terms of which
the algebraic and Hopf structures already presented in the
Serre-Chevalley basis become particularly simple. We give
furthermore the dual construction of the algebra and Hopf structure
of the functions on $SL_{qs}(1|2)$ and comment on representations of
$\ZZ_2$-graded tensor products. To further simplify calculations (by
eliminating the need to keep track of super statistics), we apply,
in the next section, the general theory of
bosonization~\cite{Majid4}, thus obtaining a standard FRT-type Hopf
algebra and its dual.
\subsection{The deformed universal enveloping algebra} \label{SFc}
We introduce the $3 \times 3$ matrices $\uL^+$, $\uL^-$ (upper
and lower triangular respectively), the elements of which, together
with 1, generate $\uqs$ and are given, in the Serre-Chevalley basis, by
\be
\begin{array}{ll}
\bigg. \uL^+_{11} = q^{H_2} s^{H_2}  \;
&\uL^-_{11} = q^{-H_2} s^{H_2}  \; \\
\bigg. \uL^+_{22} = q^{H_2-H_1} s^{H_2+H_1} \;
&\uL^-_{22} = q^{-H_2+H_1} s^{H_2+H_1} \; \\
\bigg. \uL^+_{33} = q^{-H_1} s^{H_1}  \;
&\uL^-_{33} = q^{H_1} s^{H_1}  \; \\
\bigg. \uL^+_{12} = \lambda E^+_1 q^{H_2-H_1} s^{H_2+H_1} \;
&\uL^-_{21} = -q^{-1} \lambda E^-_1 q^{-H_2+H_1} s^{H_2+H_1}  \; \\
\bigg. \uL^+_{23} = \lambda E^+_2 q^{-H_1} s^{H_1}  \;
&\uL^-_{32} = q \lambda E^-_2 q^{H_1} s^{H_1} \; \\
\bigg. \uL^+_{13} = q \lambda E^+_3 q^{-H_1} s^{H_1}  \;
&\uL^-_{31} = q \lambda E^-_3 q^{H_1} s^{H_1}  \; .
\end{array}
\label{LpmSC}
\ee
The relations\men{CommQuant},\men{SerreRelQuant} can now be
expressed compactly as follows
\be
R_{12} \uL^{\pm}_{2} \eta_{12} \uL^{\pm}_{1} =
\uL^{\pm}_{1} \eta_{12} \uL^{\pm}_{2} R_{12}, \; \; \; \; \;
R_{12} \uL^{+}_{2} \eta_{12} \uL^{-}_{1} =
\uL^{-}_{1} \eta_{12} \uL^{+}_{2} R_{12},
\label{RLeL}
\ee
where
\be
\eta_{ik,jl} = (-1)^{ik} \delta_{ij} \delta_{kl} = diag(
-1,1,-1,1,1,1,-1,1,-1)
\label{etanum}
\ee
and the $\uqs$ $R$-matrix is given by
\be
R=
\left(
\begin{array}{ccc|ccc|ccc}
-1 & 0 & 0 & 0 & 0 & 0 & 0 & 0 & 0 \\
0 & q^{-1}s & 0 & 0 & 0 & 0 & 0 & 0 & 0 \\
0 & 0 & -q^{-1}s & 0 & 0 & 0 & 0 & 0 & 0 \\
\hline
0 & -q^{-1} \lambda & 0 & q^{-1} s^{-1} & 0 & 0 & 0 & 0 & 0 \\
0 & 0 & 0 & 0 & q^{-2} & 0 & 0 & 0 & 0 \\
0 & 0 & 0 & 0 & 0 & q^{-1}s & 0 & 0 & 0 \\
\hline
0 & 0 & -q^{-1}\lambda & 0 & 0 & 0 & -q^{-1} s^{-1} & 0 & 0 \\
0 & 0 & 0 & 0 & 0 & -q^{-1}\lambda & 0 & q^{-1} s^{-1} & 0 \\
0 & 0 & 0 & 0 & 0 & 0 & 0 & 0 & -1
\end{array}
\label{Rexp}
\right)
\ee
with the matrix $\eta$ as its classical limit.
The explicit form of\men{RLeL} can be found in the appendix.
In the following we will make frequent use of the relations
\be
\eta^2 = I, \; \; \;
\eta_{12} = \eta_{21}, \; \; \;
\eta R = R \eta, \; \; \;
R_{12} \eta_{13} \eta_{23} = \eta_{23} \eta_{13} R_{12}.
\label{etaR}
\ee
Starting from the product $\uL^{\pm}_{1} \eta_{12}
\uL^{\pm}_{2} \eta_{13} \eta_{23} \uL^{\pm}_{3}$ and
using\men{RLeL}, in two different ways, to bring it to the form
$\uL^{\pm}_{3} \eta_{23}
\eta_{13} \uL^{\pm}_{2} \eta_{12} \uL^{\pm}_{1}$
(according to the sequence of transpositions $123 \ra 132 \ra
312 \ra 321$ and $123 \ra 213 \ra 231 \ra 321$ respectively), one
ensures that no cubic relations are imposed among the
$\uL^{\pm}_{ij}$'s if $R$ satisfies the (ordinary) Quantum
Yang-Baxter Equation (QYBE)
\be
R_{12} R_{13} R_{23} = R_{23} R_{13} R_{12} \; ; \label{QYBE}
\ee
it is easily verified that the $R$ given by\men{Rexp} indeed
satisfies this relation. In terms of the matrix $\hR$, given by
$\hR_{ik,jl} = R_{ki,jl}$,\men{QYBE} reads
\be
\hRot \hRtth \hRot = \hRtth \hRot \hRtth \label{QYBEh}.
\ee
The $\hR$ matrix is non-symmetric in our case and it satisfies the
quadratic characteristic equation
\be
\hR^{2} + q^{-1} \lambda \hR - \qit =0 \; ,
\label{Rhcheqn}
\ee
which implies the eigenvalues $-1$ and $\qit$. With the fermionic
degree of $\uL^{\pm}_{ij}$ being given by
\be
\deg \uL^{\pm}_{ij} = i+j \quad (\mod 2) \; ,
\label{degLpm}
\ee
we can express the superstatistics in matrix form as follows
\be
\uL^{\pm}_{1} \eta_{12} {\uL^{\pm}_{2}}' \eta_{12} =
\eta_{12} {\uL^{\pm}_{2}}' \eta_{12} \uL^{\pm}_{1} \; ,
\; \; \; \; \; \;
\uL^{+}_{1} \eta_{12} {\uL^{-}_{2}}' \eta_{12} =
\eta_{12} {\uL^{-}_{2}}' \eta_{12} \uL^{+}_{1} \; ,
\label{LLstat}
\ee
where $\uL^{\pm}$, ${\uL^{\pm}}'$ denote two copies of the generators
living in different spaces in the tensor product $\uqs
\uot \uqs$. The Hopf structure given in section\mens{Quantum I} is
now compactly encoded in the relations
\be
\uDelta (\uLpm ) = \uLpm \dot{\uot} \uLpm \; ,
\; \; \; \; \;
\us (\uLpm) = (\uLpm)^{-1},
\; \; \; \; \; \;
\ueps (\uLpm) = I \; ,
\ee
where $\dot{\uot}$ stands for $\ZZ_{2}$-graded tensor product and
matrix multiplication. Notice finally the superdeterminant
relations (obtained by inspection of\men{LpmSC})
\be
\uL^{\pm}_{11} (\uL^{\pm}_{22})^{-1} \uL^{\pm}_{33} = 1 \; .
\label{sdet}
\ee
One can combine all the generators in a single matrix $\uY$,
given by~\cite{FRT}
\be
\uY = \uLp \us (\uLm) \; ,
\label{Ydef}
\ee
the elements of which are easily seen to satisfy the commutation
relations
\be
R_{21} \uY_{1} R_{12} \uY_{2} = \uY_{2} R_{21} \uY_{1} R_{12} \; ,
\label{RYRY}
\ee
with $\deg \uY_{ij} = i+j \quad (\mod 2)$.
\subsection{The dual superalgebra of functions}
\label{dsf}
\subsubsection{The construction of the dual}

\indent

The algebra $\Aqs = {\rm Fun}(SL_{qs}(1|2))$, dual to $\uqs$,
 is generated
by 1 and the elements of the ($3 \times 3$, in our case) supermatrix
$\uA$ -- their inner product with $\uLpm$ is given by
\be
\bip{\uLp_{1}}{\uA_{2}} = \eta_{12} R_{21} \; ,
\; \; \; \; \;
\bip{\uLm_{1}}{\uA_{2}} = \eta_{12} R^{-1}_{12} \; ,
\; \; \; \; \;
\bip{\uU}{\uA} = I \; .
\label{LpmAip}
\ee
The $\ZZ_{2}$-graded Hopf structure is
\be
\uDelta (\uA) = \uA \dot{\uot} \uA \; ,
\; \; \; \; \; \;
\uS (\uA) = \uA^{-1} \; ,
\; \; \; \; \; \;
\ueps (\uA) = I \; ,
\label{AHs}
\ee
with the coproduct being related to the product in the dual via
($x$, $y$ $\in$ $\uqs$)
\be
\bip{xy}{\uA_{ij}} = \bip{x}{\uA_{im}} \bip{y}{\uA_{mj}} \; .
\label{UApcop}
\ee
The reader will notice here that although going from the lhs to the
rhs of\men{UApcop} involves (typographicaly) a transposition of $y$
and $\uA_{im}$, no corresponding sign factor was included on the
rhs. The advantage of this convention is that it supplies us with a
matrix representation of $\uqs$ (with the usual matrix
multiplication), a feature that we find especially appealing in
view of the applications of section\mens{Hamil}. On the
other hand, viewing this quantum supergroup as a particular example
of a braided group (see, for example,~\cite{Majid4} and references
therein), one would probably want to attempt a formulation that is
compatible with the diagrammatics of category theory. We give, for
completeness, such a formulation at the end of this section.

We check now, in order to illustrate the formalism,
whether the first of\men{AHs} is compatible with the
``$RLL$'' commutation relations, Eq.\men{RLeL}. The inner
product of the lhs of\men{RLeL} with $\uA$ is given by
\bae
\bip{R_{23} \uLp_{3} \eta_{23} \uL^{+}_{2}}{\uA_{1}} \fe
R_{23} \bip{\uL^{+}_{3}}{\uA_{1}} \eta_{23} \bip{\uL^{+}_{2}}{\uA_{1}}
\ff
 \fe R_{23} \eta_{13} R_{13} \eta_{23} \eta_{12} R_{12} \ff
 \fe R_{23} R_{13} \eta_{13} \eta_{23} R_{12} \eta_{12} \ff
 \fe R_{23} R_{13} R_{12} \eta_{12} \eta_{13} \eta_{23} \; , \nn
\eae
while for the right hand side we get
\bae
\bip{\uL^{+}_{2} \eta_{23} \uL^{+}_{3} R_{23}}{\uA_{1}}
 \fe \bip{\uLp_{2}}{\uA_{1}} \eta_{23} \bip{\uLp_{3}}{\uA_{1}} R_{23}
\ff
 \fe \eta_{12} R_{12} \eta_{23} \eta_{13} R_{13} R_{23} \ff
 \fe R_{12} \eta_{12} \eta_{23} R_{13} \eta_{13} R_{23} \ff
 \fe R_{12} R_{13} \eta_{12} \eta_{13} \eta_{23} R_{23} \ff
 \fe R_{12} R_{13} R_{23} \eta_{12} \eta_{13} \eta_{23} \nn
\eae
and, by invoking the QYBE for $R$, one verifies that $\uA$ is a
representation.
 The degree of its elements is given by
\be
\deg \uA_{ij} = (-1)^{i+j} \quad (\mod 2) \; ;
\label{degA}
\ee
(this is the standard fermionic basis format - the
block format is standard in the distinguished basis).
Notice that the requirement that
the quantum superdeterminants
of Eq.\men{sdet} be represented by the unit matrix fixes the
normalization of the $R$-matrix.
\subsubsection{Tensor product representations}
We address next the question of tensor product representations.
It is clear that the matrix $\uA_{1} \uot \uA_{2}$ does not provide
a representation for $\uqs \uot \uqs$ since, for example, the elements
$\uLp \uot 1$ and $1 \uot \uLp$ do not commute (see
Eq.\men{LLstat}), while their inner products with $\uA_{1} \uot
\uA_{2}$ necessarily do. An alternative statement of this fact can
be made via the naturally induced coproduct $\uDelta^{ \{ 2 \} }$
in $\Aqs \uot \Aqs$
\be
\uDelta^{\{ 2 \} } (a \uot b) =
(-1)^{\deg b_{(\underline{1})} \, \deg a_{(\underline{2})}}
\; (a_{(\underline{1})} \uot b_{(\underline{1})})
\uot (a_{(\underline{2})} \uot b_{(\underline{2})})
\label{DDelta}
\ee
where (Sweedler notation)
\be
\uDelta(a) \equiv \sum_{i} a_{(\underline{1})}^{i} \uot a_{(\underline{2})}^{i}
\equiv a_{(\underline{1})} \uot a_{(\underline{2})} \; . \nn
\ee
The failure of $\uA_{1} \uot \uA_{2}$ to provide a representation
for $\uqs \uot \uqs$ can be traced to the fact that
\be
\uDelta^{ \{ 2 \} }(\uA_{1} \uot \uA_{2}) \neq (\uA_{1} \uot \uA_{2})
\uot (\uA_{1} \uot \uA_{2}) \; . \nn
\ee
We find though that the matrix
\be
\uA_{12} \equiv \uA_{1} \uot \eta_{12} \uA_{2} \eta_{12} \label{A12def}
\ee
satisfies
\be
\uDelta^{ \{ 2 \} } (\uA_{12}) = \uA_{12} \uot \uA_{12} \; . \label{A12cop}
\ee
It is clear that $\uA_{12}$ is a representation for $\uqs \uot \uU$ and
$\uU \uot \uqs$ (the latter since $\eta^{2} = I$). We check whether
it represents the ``cross'' commutation relations, Eq.\men{LLstat}.
For the inner product of the lhs of the first of\men{LLstat} with
$\uA_{34}$ we get
\bae
\bip{(\uLp_{1} \uot \uU) \eta_{12} ( \uU \uot \uLp_{2})
\eta_{12}}{\uA_{34}}
 \fe \bip{\uLp_{1} \uot \uU}{\uA_{34}} \eta_{12} \bip{\uU \uot
\uLp_{2}}{\uA_{34}} \eta_{12} \ff
 \fe R_{31} \eta_{34}^{2} \eta_{12} \eta_{34} R_{42} \eta_{34}
\eta_{12} \ff
 \fe R_{31} \eta_{12} \eta_{23}^{2} \eta_{34} R_{42} \eta_{34}
\eta_{12} \ff
 \fe \eta_{12} \eta_{23} R_{31} R_{42} \eta_{23} \eta_{12} \nn
\eae
while the rhs gives
\bae
\bip{\eta_{12} ( \uU \uot \uLp_{2}) \eta_{12} (\uLp_{1} \uot
\uU)}{\uA_{34}}
 \fe \eta_{12} \bip{\uU \uot \uLp_{2}}{\uA_{34}} \eta_{12}
\bip{\uLp_{1} \uot \uU}{\uA_{34}} \ff
 \fe \eta_{12} \eta_{34} R_{42} \eta_{34} \eta_{12} R_{31} \ff
 \fe \eta_{12} \eta_{34} R_{42} \eta_{34} \eta_{23} \eta_{23}
\eta_{12} R_{31} \ff
 \fe \eta_{12} \eta_{34} \eta_{34} \eta_{23} R_{42} R_{31}
\eta_{23} \eta_{12} \ff
 \fe \eta_{12} \eta_{23} R_{42} R_{31} \eta_{23} \eta_{12} \; ; \nn
\eae
we conclude that $\uA_{34}$ provides a representation for $\uqs \uot \uqs$.
The extension to higher tensor products gives the representation
matrix $\uA_{1 \ldots L}$ defined by
\be
\uA_{1 \ldots L} = \bigotimes\limits_{\overline{\phantom{.}
\atop \scriptstyle{k=1}}}^{L} \uA_{(k)}
\label{A1Ldef}
\ee
where
\be
\uA_{(k)} \equiv \eta_{1k} \ldots \eta_{k-1,k} \uA_{k} \eta_{k-1,k}
\ldots \eta_{1k}
\label{Akdef}
\ee
which satisfies
\be
\uDelta^{ \{ L \} } (\uA_{1 \ldots L}) = \uA_{1 \ldots L} \uot \uA_{1
\ldots L} \; ,
\label{DA1L}
\ee
$\uDelta^{ \{ L \} }$ being the naturally induced coproduct
in $\Aqs^{\uot
L}$.

A number of consistency checks are in order at this point. We
verify, as an example, the compatibility of the $\uA-\uA$ commutation
relations with the $\uA-\uA'$ superstatistics (given by an equation
analogous to\men{LLstat})
\bae
(R_{12} \uA_{1} \eta_{12} \uA_{2} - \uA_{2} \eta_{12} \uA_{1}
R_{12} ) \eta_{13} \eta_{23}  \uA_{3}'
 \fe R_{12} \uA_{1} \eta_{12} \eta_{13} \eta_{23} \uA_{3}'
\eta_{23} \uA_{2} \eta_{23} \ff
 & & \; \;  - \uA_{2} \eta_{12} \eta_{23} \eta_{13}
\uA_{3}' \eta_{13} \uA_{1} \eta_{13} R_{12} \ff
 \fe R_{12} \eta_{23} \eta_{13} \uA_{3}' \eta_{13} \uA_{1}
\eta_{13} \eta_{12} \eta_{23} \uA_{2} \eta_{23} \ff
  & & \; \;  - \eta_{13}
\eta_{23} \uA_{3}' \eta_{23} \uA_{2} \eta_{23} \eta_{12} \eta_{13} \uA_{1}
\eta_{13}
R_{12} \ff
 \fe \eta_{23} \eta_{13} \uA_{3}' R_{12} \eta_{13} \eta_{23}
\uA_{1} \eta_{12} \uA_{2} \eta_{13} \eta_{23} \ff
 & & \; \;  - \eta_{13} \eta_{23} \uA_{3}'
\eta_{23} \eta_{13} \uA_{2} \eta_{12}
\uA_{1} R_{12} \eta_{13} \eta_{23} \ff
 \fe \eta_{13} \eta_{23} \uA_{3}' \eta_{13} \eta_{23}
( R_{12} \uA_{1} \eta_{12} \uA_{2} -
\uA_{2} \eta_{12} \uA_{1} R_{12}) \eta_{13} \eta_{23} \; . \nn
\eae
One can show, along the same lines, the consistency of the
entire scheme presented above.
\subsubsection{$\ZZ_{2}$-graded inner product}
As promised earlier in this section, we present now a version of
the construction of $\Aqs$ that uses the $\ZZ_{2}$-graded inner
product $\tip{\cdot}{\cdot}$ defined by
\be
\tip{\uLp_{1}}{\eta_{12} \uA_{2}} = R_{21} \; ,
\; \; \; \; \; \;
\tip{\uLm_{1}}{\eta_{12} \uA_{2}} = R_{12}^{-1} \; ,
\; \; \; \; \; \;
\tip{1}{\uA} = I
\label{UApcop2}
\ee
and the duality relations ($x$, $y \in \uqs$, $a$, $b \in \Aqs$,
with $y$, $a$ of homogeneous degree)
\bae
\tip{xy}{\uA_{ij}} \fe \tip{x \uot y}{\uA_{im} \uot \uA_{mj}} \ff
 & \equiv & (-1)^{\deg{\uA_{im}} \cdot \deg{y}} \tip{x}{\uA_{im}}
 \tip{y}{\uA_{mj}} \; , \\
\tip{\uLpm_{ij}}{ab} \fe \tip{\uLpm_{im} \uot \uLpm_{mj}}{a \uot b}
\ff
 & \equiv & (-1)^{\deg{\uLpm_{mj}} \cdot \deg{a}}
\tip{\uLpm_{im}}{a} \tip{\uLpm_{mj}}{b} \; .
\eae
$\Aqs$ remains unchanged as a Hopf algebra. One can go further and
give the natural semidirect (graded) product commutation relations
between $\uLpm$ and $\uA$
\bae
\uLp_{1} \eta_{12} \uA_{2} \fe \eta_{12} \uA_{2} R_{21} \uLp_{1}
\eta_{12} \label{LpAcr} \\
\uLm_{1} \eta_{12} \uA_{2} \fe \eta_{12} \uA_{2} R^{-1}_{12} \uLm_{1}
\eta_{12} \; .
\eae
\sect{Bosonization} \label{Boson}
\subsection{A bosonised basis}

\indent

We start with the following observation: the element $g = e^{i \pi
(H_{1} + H_{2})}$ of $\uqs$ satisfies
\be
g \, x = (-1)^{\deg x} x \, g \label{gx}
\ee
where $x$ is any element (of homogeneous degree) of $\uqs$. Notice
also that $\uDelta(g) = g \uot g$ and
\be
\bip{g}{\uA} = \left( \begin{array}{ccc}
                  -1 & 0 & 0 \\
                  0 & 1 & 0 \\
                  0 & 0 & -1
                \end{array} \right)
\ee
so that one can consistently set $g^{2} = 1$ as an operator
relation in the algebra. Under these conditions, one may apply the
general construction of bosonization of a super-Hopf
algebra~\cite{Majid4} to obtain an ordinary (\ie \,
non-$\ZZ_{2}$-graded) coproduct for the $\uLpm_{ij}$ 's
\be
\Delta (\uLpm_{ij}) = \sum_{m} \uLpm_{im} g^{m+j} \otimes \uLpm_{mj}
\; , \label{uLcop}
\ee
along with a matching counit and antipode. The entire bosonic Hopf
structure is most compactly expressed in terms of the matrices
$\Lpm$ given by
\be
\Lpm = \uLpm G \; ,
\; \; \; \; \; \;
G \equiv \diag (g \; , 1 \; ,g) \; . \label{LuL}
\ee
Then\men{uLcop} and the rest of the Hopf structure become
\be
\Delta(\Lpm) = \Lpm \dot{\ot} \Lpm \; ,
\; \; \; \; \; \;
S(\Lpm) = {\Lpm}^{-1} \; ,
\; \; \; \; \; \;
\epsilon(\Lpm) = I \; .
\label{Lhopf}
\ee
Starting from\men{RLeL}, one finds for the $L-L$ commutation
relations the standard FRT expressions
\be
R_{12} \Lpm_{2} \Lpm_{1} = \Lpm_{1} \Lpm_{2} R_{12} \; ,
\; \; \; \; \; \;
R_{12} \Lp_{2} \Lm_{1} = \Lm_{1} \Lp_{2} R_{12} \; .
\label{RLL}
\ee
Notice that\men{RLL} just expresses the original $\uqs$
algebra in a still different basis while the Hopf structure
of\men{Lhopf} is not equivalent to the ones presented in the
previous two sections. The explicit form of\men{RLL} can be found
in the appendix. For the square of the antipode we find, by direct
calculation,
\be
S^{2}(\Lpm) = g \Lpm g^{-1} = D \Lpm D^{-1}
\label{S2Lpm}
\ee
where $D = \diag (-1 \; , 1 \; , -1)$. Since $D^2 = I$, we have
$S^4 = \id$. The $D$ matrix further
satisfies
\be
\Tr_{2} (D_{2} \hR_{12}) = I_{1} \; ,
\; \; \; \; \; \;
\Tr_{1}(D^{-1}_{1} \hR_{12}^{-1}) = I_{2} \; ,
\; \; \; \; \; \;
D_{1} (R^{T_{2}})^{-1}_{12} =  (R^{-1})^{T_{2}}_{12} D_{1} \; ,
\label{Dident}
\ee
where $T_{2}$ denotes transposition in the second matrix space.
Using the characteristic equation for $\hR$, one can derive
from\men{Dident} two more useful identities
\be
\Tr_{2} (D_{2} \hR^{-1}_{12}) = I_{1} \; ,
\; \; \; \; \; \;
\Tr_{1}(D_{1}^{-1} \hR_{12}) = I_{2} \; .
\label{Dident2}
\ee
Finally, the $Y$ matrix, defined by $Y=\Lp S(\Lm)$, is found to be
equal to $\uY$.

The construction of the dual proceeds as in~\cite{FRT} -- the
formulas for the commutation relations and the Hopf structure can
be obtained from those of the previous section by setting $\eta =
I$ and omitting the underlining of symbols. The inner product
relations are
\be
\ip{\Lp_{1}}{A_{2}} = R_{21} \; ,
\; \; \; \; \; \;
\ip{\Lm_{1}}{A_{2}} = R_{12}^{-1} \; ,
\; \; \; \; \; \;
\ip{1}{A} = I \; .
\label{ipLA}
\ee
Notice the use of a (still) different symbol for the inner product. The
ones of the previous section relate, via
equations\men{UApcop},\men{UApcop2},
products in $\uqs$ with graded coproducts in $\Aqs$ while the one
appearing in\men{ipLA} conforms to the (standard) Hopf algebra
duality requirement
\be
\ip{xy}{a} = \ip{x}{a_{(1)}} \ip{y}{a_{(2)}} \; ,
\; \; \; \; \; \;
\Delta{a} \equiv a_{(1)} \ot a_{(2)} \; . \nn
\ee
\subsection{The centre of $\uqs$}
We now turn our attention to the centre of $\uqs$. We consider the
set of Casimirs given by
\be
{c^{(k)}} = \lambda^{-k} \Tr (D^{-1} (I-Y)^{k})
\; \; \; \; \; \;
k \in \ZZ
\label{ckpdef}
\ee
(these are linearly related to the Casimirs introduced
in~\cite{FRT}). The proof of their centrality relies on the third
of\men{Dident}.
With the representation of $Y$ being given by
\be
\ip{Y_{1}}{A_{2}} = \hR_{12}^{2}
\label{Yrep}
\ee
and the characteristic equation\men{Rhcheqn} allowing us to write
\be
(I - \hR^{2})^{k} = \frac{(1-q^{-4})^{k}}{1+ q^{-2}}
(\hR + I) \; ,
\label{Rk}
\ee
we find for the values of the casimirs in the fundamental
representation
\bae
\ip{c^{(k)}}{A_{2}} \fe \Tr_{1} (D_{1}^{-1}
\ip{(I_{1} - Y_{1}))^{k}}{A_{2}}) \ff
 \fe \Tr_{1} (D_{1}^{-1} (I_{12} - \hR_{12}^{2})^{k}) \ff
 \fe \frac{(1-q^{-4})^{k}}{1+ q^{-2}} \Tr_{1} (D_{1}^{-1}(\hR_{12}
 + I_{12}))
\ff
\Rightarrow \ip{c^{(k)}}{A_{2}} \fe 0  \; .
\eae
The construction of the Hamiltonians of the next section requires
the evaluation of the representation of the coproduct of central
elements in $\uqs$. For the above set of Casimirs (and $k \geq 1$)
 we compute
\bae
c^{(k)}_{23} & \equiv & \ip{\Delta(c^{(k)})}{A_{2} \ot A_{3}} \ff
 \fe \Tr_{1} D_{1}^{-1} \lambda^{-k}
\ip{(\Delta(I_{1}-Y_{1}))^{k}}{A_{2} \ot A_{3}} \ff
 \fe \Tr_{1} D_{1}^{-1} X_{123}^{k}
\label{Dckrep}
\eae
where
\bae
X_{123} & \equiv & \lambda^{-1} \ip{\Delta(I_{1} - Y_{1})}{A_{2} \ot
A_{3}} \ff
 \fe q^{-1} \hR_{23} \hR_{12} \hR_{23} + q^{-3} \hR_{12} + q^{-1} +
q^{-3} \; .
\label{X123def}
\eae
To extract more information about $c^{(k)}_{23}$, we
assume that, for some $k$, $X_{123}^{k}$ is of the form
\be
X_{123}^{k} = a_{k} \hR_{23} \hR_{12} \hR_{23}
+ b_{k} (\hR_{23} \hR_{12} + \hR_{12} \hR_{23})
 +c_{k} \hR_{12} + d_{k} \hR_{23} + f_{k} I
\label{Xkform}
\ee
and, furthermore, that for that same $k$, the relation
\be
c_{k} = q^{-2} a_{k} + (1-q^{-2}) b_{k}
\label{ckrel}
\ee
holds (both assumptions are evidently valid for $k=1$). Then,
multiplying\men{Xkform} by $X_{123}$ and using the characteristic
equation for $\hR$, we find the recursion relations
\bae
a_{k+1} \fe (q^{-7} - q^{-5} +2 q^{-3} ) a_{k}
 + (2 q^{-5} -3 q^{-3} + 2q^{-1} ) b_{k} + (q^{-3} - q^{-1} ) c_{k}
\ff
 & & \; \; + (q^{-3} - q^{-1}) d_{k} + q^{-1} f_{k} \ff
b_{k+1} \fe (q^{-7} - q^{-5} + q^{-3}) a_{k} +
(2 q^{-5} - q^{-3} + q^{-1}) b_{k} + q^{-3} c_{k} + q^{-3} d_{k} \ff
c_{k+1} \fe (q^{-7} - q^{-5}) a_{k} + q^{-5} b_{k} +
(q^{-5} + q^{-1}) c_{k}
+ q^{-3} f_{k} \ff
d_{k+1} \fe (q^{-7} - q^{-5}) a_{k} + 2 q^{-5} b_{k} +
(q^{-3} + q^{-1}) d_{k} \ff
f_{k+1} \fe q^{-7} a_{k} + q^{-5} c_{k} + (q^{-3} + q^{-1}) f_{k}
\label{recrel}
\eae
from which it easily follows, by induction, that both\men{Xkform}
and\men{ckrel} are valid for all $k \geq 1$. Taking the quantum
trace on both sides of\men{Xkform} we find, with the help of the
second of\men{Dident}
\be
c^{(k)}_{23} = (-q^{-1} \lambda a_{k} + 2b_{k} - d_{k}) \hR_{23}
 + (q^{-2} a_{k} + c_{k} - f_{k}) I \; .
\label{ck23one}
\ee
We assume that, for some $k$, the above expression is proportional
to $\hR_{23} + I$ (true for $k=1$). This gives us the relation
\be
f_{k} = (1+q^{-2}) a_{k} - (1+q^{-2}) b_{k} + d_{k}
\label{fkrel}
\ee
which is easily established, by induction, for all $k \geq 1$.
Implementing\men{ckrel},\men{fkrel} in\men{recrel} we are left with
a recursion relation for the sequences $a_{k}$, $b_{k}$, $d_{k}$ --
solving it and substituting in\men{ck23one} we finally find
\be
c^{(k)}_{23} = \alpha_{k}(\hR+I)_{23} \; ,
\; \; \; \; \; \; \; \; \;
\alpha_{k} = \frac{q^{-4k+1}}{[ 2 ]_{q} [ 3 ]_{q}}
([ 4 ]_{q}^{k} - q^{3k} [ 2 ]_{q}^{2})
 \; . \label{ck23two}
\ee

\sect{Hamiltonians \label{Hamil}}

\subsection{General scheme}

\indent

In this section, we apply the results concerning the Casimir operators
to construct $3^L$-state quantum chain Hamiltonians with nearest
neighbour interaction that are
$sl(1|2)$ (resp. $\cU_{qs}(sl(1|2))$) invariant.

Define a $L$-site Hamiltonian $\cH^{(1 \cdots L)}_p$ by
\be
\cH^{(1 \cdots L)}_p = \sum_{j=1}^{L-1} \cH^{(j,j+1)}_p
\ee
where
\be
\cH^{(j,j+1)}_p = 1 \otimes \cdots \otimes \cH_p \otimes \cdots \otimes 1
\ee
with the two-site Hamiltonian $\cH_p$ in position $(j,j+1)$.
\\
The two-site Hamiltonian $\cH_p$ itself is defined by
\be
(\cH_p)_{12} = \bip{\uDelta(\cC_p)}{\uA_{12}}
\ee
for any Casimir operator $\cC_p$ of $\cU$ (resp. $\uqs$).
\\
Then the $L$-site Hamiltonian $\cH^{(1 \cdots L)}_p$ is
$\cU$-invariant (resp. $\uqs$-invariant) in the sense that
\be
\left[\cH^{(1 \cdots L)}_p , X_{1 \cdots L} \right] = 0
\ee
where $X_{1 \cdots L} $ is the evaluation in the $L$-fold
tensor product of the
fundamental representation of $X\in \cU$ (resp. $\uqs$) \ie
\be
X_{1 \cdots L} \equiv \bip {\uDelta^{(L-1)} (X)}{ \uA_{1 \cdots L}} \;.
\ee

\subsection{The classical case}

\indent

The fundamental representation of $sl(1|2)$ is given by
\be
\begin{array}{lll}
\pi(H_1) = e_{11} + e_{22}\;, \quad & \pi(H_2) = -e_{22} -
e_{33}\;, \quad &
\\
\bigg. \pi(E^+_1) = e_{21}\;,&
\pi(E^+_2) = e_{32}\;,&
\pi(E^+_3) = e_{31}\;,
\\
\bigg. \pi(E^-_1) = e_{12}\;,&
\pi(E^-_2) = -e_{23}\;,&
\pi(E^-_3) = -e_{13}\;,
\end{array}
\label{FundRepClass}
\ee
where $e_{ij}$ are $3 \times 3$ elementary matrices with entry 1 in
$i$th row and $j$th column and 0 elsewhere.
\\
By applying the general scheme above to the classical Casimir
operators $\cC_p^{cl}$ computed in Section \ref{Class}, one finds that
the corresponding Hamiltonians $\cH_p^{cl}$ are of the form
$\cH_p^{cl} = \alpha_p \cH^{cl} + \beta_p I$ where $I$ is the identity and
$\cH^{cl}$ is given by
\be
\cH^{cl} =
\left(
\begin{array}{ccc|ccc|ccc}
0 & 0 & 0 & 0 & 0 & 0 & 0 & 0 & 0 \\
0 & 1 & 0 & 1 & 0 & 0 & 0 & 0 & 0 \\
0 & 0 & 1 & 0 & 0 & 0 & -1 & 0 & 0 \\
\hline
0 & 1 & 0 & 1 & 0 & 0 & 0 & 0 & 0 \\
0 & 0 & 0 & 0 & 2 & 0 & 0 & 0 & 0 \\
0 & 0 & 0 & 0 & 0 & 1 & 0 & 1 & 0 \\
\hline
0 & 0 & -1 & 0 & 0 & 0 & 1 & 0 & 0 \\
0 & 0 & 0 & 0 & 0 & 1 & 0 & 1 & 0 \\
0 & 0 & 0 & 0 & 0 & 0 & 0 & 0 & 0
\end{array}
\right) \;.
\label{HamilClass}
\ee

To get a physical interpretation of $\cH^{cl}$, one can use the
realization of $sl(1|2)$ in terms of fermionic creation and
annihilation operators
$c_{\uparrow i}^{\pm}$ and $c_{\downarrow i}^{\pm}$ such that
\be
\{c_{\alpha i}^-,c_{\beta j}^+ \} = \delta_{\alpha \beta} \delta_{ij}
\ee
where $\alpha, \beta = \uparrow$ or $\downarrow$ and $i,j$ are the site
indices.
\\
Let us define the spin operators $\sigma_j^+, \sigma_j^-, \sigma_j^0$ by
\be
\sigma_j^+ = c_{\uparrow j}^+ c_{\downarrow j}^- \;, \ \ \
\sigma_j^- = c_{\downarrow j}^+ c_{\uparrow j}^- \;, \ \ \
\sigma_j^0 = c_{\uparrow j}^+ c_{\uparrow j}^-
- c_{\downarrow j}^+ c_{\downarrow j}^- \;,
\label{SigmaOp}
\ee
and the number operators $n_{\uparrow j}$, $n_{\downarrow j}$, $n_j$,
$n_{\emptyset j}$ by
\be
n_{\uparrow j} = c_{\uparrow j}^+ c_{\uparrow j}^- \;, \ \ \
n_{\downarrow j} = c_{\downarrow j}^+ c_{\downarrow j}^- \;, \ \ \
n_j = n_{\uparrow j} + n_{\downarrow j} \;, \ \ \
n_{\emptyset j}= 1 - n_j \;.
\label{NumbOp}
\ee
Then the realization of the superalgebra $sl(1|2)$ is given by
\be
\begin{array}{lll}
e_{11} = n_{\uparrow} \;, & e_{22} = n_{\emptyset}\;,
& e_{33} = n_{\downarrow} \;,
\\
\bigg. e_{12} = c_{\uparrow}^+ (1-n_{\downarrow}) \;, \quad &
e_{23} = c_{\downarrow}^- (1-n_{\uparrow}) \;, \quad &
e_{13} = \sigma^+ \;,
\\
\bigg. e_{21} = c_{\uparrow}^- (1-n_{\downarrow}) \;,&
e_{32} = c_{\downarrow}^+ (1-n_{\uparrow}) \;,&
e_{31} = \sigma^- \;.
\end{array}
\label{Realiz}
\ee
It follows that the two-site Hamiltonian ${\cH^{cl}}^{(j,j+1)}$
can be expressed as
\bea
{\cH^{cl}}^{(j,j+1)} &=& + \ c_{\uparrow j}^+ (1 - n_{\downarrow j})
       c_{\uparrow {j+1}}^- (1 - n_{\downarrow {j+1}})
    -  c_{\uparrow j}^- (1 - n_{\downarrow j})
       c_{\uparrow{j+1}}^+ (1 - n_{\downarrow {j+1}}) \nonumber \\
& & +\ c_{\downarrow j}^+ (1 - n_{\uparrow j})
       c_{\downarrow {j+1}}^- (1 - n_{\uparrow {j+1}})
    -  c_{\downarrow j}^- (1 - n_{\uparrow j})
       c_{\downarrow{j+1}}^+ (1 - n_{\uparrow {j+1}}) \nonumber \\
& & -\ \sigma_j^+ \sigma_{j+1}^- - \sigma_j^- \sigma_{j+1}^+ \nonumber \\
& & +\ (1 - n_{\downarrow j})(1 - n_{\uparrow {j+1}})
    + (1 - n_{\uparrow j})(1 - n_{\downarrow {j+1}}) \;.
\label{HamiltJClass}
\ena
which is precisely the Hamiltonian that describes the $t-J$ model at
the supersymmetric point. This model was found to be invariant
under the $sl(1|2)$ superalgebra \cite{WFS}.

\subsection{The deformed case}

\subsubsection{$\uqs$-invariant Hamiltonians from the fermionic basis}

\indent

The fundamental representation of $\uqs$ remains, as usual, undeformed
and is given by (\ref{FundRepClass}).
Using the FRT formalism and the set of Casimir operators $c^{(k)}$ of
Sect. \ref{Boson}, one gets $\cH_{12}=c^{(k)}_{12}$, see Eq. (\ref{ck23two}).
One can now compute the Hamiltonians $\cH_p$ associated to the
Casimir operators $\cC_p$. The result is
\be
\cH_p = -q^{3-6p} (q-q^{-1})^2 \cH_{ferm} \;,
\ee
where
\be
\cH_{ferm} =
\left(
\begin{array}{ccc|ccc|ccc}
0 & 0 & 0 & 0 & 0 & 0 & 0 & 0 & 0 \\
0 & q^{-1} & 0 & s^{-1} & 0 & 0 & 0 & 0 & 0 \\
0 & 0 & q^{-1} & 0 & 0 & 0 & -s^{-1} & 0 & 0 \\
\hline
0 & s & 0 & q & 0 & 0 & 0 & 0 & 0 \\
0 & 0 & 0 & 0 & q + q^{-1} & 0 & 0 & 0 & 0 \\
0 & 0 & 0 & 0 & 0 & q^{-1} & 0 & s^{-1} & 0 \\
\hline
0 & 0 & -s & 0 & 0 & 0 & q & 0 & 0 \\
0 & 0 & 0 & 0 & 0 & s & 0 & q & 0 \\
0 & 0 & 0 & 0 & 0 & 0 & 0 & 0 & 0
\end{array}
\right) \;.
\label{HamilQuant}
\ee
This Hamiltonian is again proportional to $\hR+I$.
\\
As before, the two-site Hamiltonian $\cH^{(j,j+1)}$
can be expressed in terms of the fermionic creation and
annihilation operators
$c_{\uparrow i}^{\pm}$ and $c_{\downarrow i}^{\pm}$. One finds
\bea
\cH_{ferm}^{(j,j+1)} &=& + \ s^{-1} c_{\uparrow j}^+ (1 - n_{\downarrow j})
       c_{\uparrow {j+1}}^- (1 - n_{\downarrow {j+1}})
    -  s\ c_{\uparrow j}^- (1 - n_{\downarrow j})
       c_{\uparrow{j+1}}^+ (1 - n_{\downarrow {j+1}}) \nonumber \\
& & +\ s\ c_{\downarrow j}^+ (1 - n_{\uparrow j})
       c_{\downarrow {j+1}}^- (1 - n_{\uparrow {j+1}})
    -  s^{-1} c_{\downarrow j}^- (1 - n_{\uparrow j})
       c_{\downarrow{j+1}}^+ (1 - n_{\uparrow {j+1}}) \nonumber \\
& & -\ s^{-1} \sigma_j^+ \sigma_{j+1}^-
    - s\ \sigma_j^- \sigma_{j+1}^+ \nonumber \\
& & +\ q^{-1}(1 - n_{\downarrow j})(1 - n_{\uparrow {j+1}})
    +  q (1 - n_{\uparrow j})(1 - n_{\downarrow {j+1}}) \;.
\label{HamiltJQuant}
\ena

\subsubsection{A four parametric Hamiltonian}

\medskip
One can remark that one could have started in Sect. \ref{quantum II}
with the following more general four-parameter $R$-matrix
\be
R=
\left(
\begin{array}{ccc|ccc|ccc}
-1 & 0 & 0 & 0 & 0 & 0 & 0 & 0 & 0 \\
0 & q^{-1}q_{12} & 0 & 0 & 0 & 0 & 0 & 0 & 0 \\
0 & 0 & -q^{-1}q_{13} & 0 & 0 & 0 & 0 & 0 & 0 \\
\hline
0 & -q^{-1} \lambda & 0 & q^{-1} q_{12}^{-1} & 0 & 0 & 0 & 0 & 0 \\
0 & 0 & 0 & 0 & q^{-2} & 0 & 0 & 0 & 0 \\
0 & 0 & 0 & 0 & 0 & q^{-1}q_{23} & 0 & 0 & 0 \\
\hline
0 & 0 & -q^{-1}\lambda & 0 & 0 & 0 & -q^{-1} q_{13}^{-1} & 0 & 0 \\
0 & 0 & 0 & 0 & 0 & -q^{-1}\lambda & 0 & q^{-1} q_{23}^{-1} & 0 \\
0 & 0 & 0 & 0 & 0 & 0 & 0 & 0 & -1
\end{array}
\label{Rqij}
\right) \;.
\ee
Using the fundamental representation and the Casimir operators
of the algebra defined by the relations (\ref{RLeL}), with this
$R$-matrix (\ref{Rqij}), one gets the following quantum chain
Hamiltonian
\be
\cH =
\left(
\begin{array}{ccc|ccc|ccc}
0 & 0 & 0 & 0 & 0 & 0 & 0 & 0 & 0 \\
0 & q^{-1} & 0 & q_{12}^{-1} & 0 & 0 & 0 & 0 & 0 \\
0 & 0 & q^{-1} & 0 & 0 & 0 & q_{13}^{-1} & 0 & 0 \\
\hline
0 & q_{12} & 0 & q & 0 & 0 & 0 & 0 & 0 \\
0 & 0 & 0 & 0 & q + q^{-1} & 0 & 0 & 0 & 0 \\
0 & 0 & 0 & 0 & 0 & q^{-1} & 0 & q_{23}^{-1} & 0 \\
\hline
0 & 0 & q_{13} & 0 & 0 & 0 & q & 0 & 0 \\
0 & 0 & 0 & 0 & 0 & q_{23} & 0 & q & 0 \\
0 & 0 & 0 & 0 & 0 & 0 & 0 & 0 & 0
\end{array}
\right) \; ,
\label{Hamilqij}
\ee
which reads, in terms of fermionic operators,
\bea
\cH^{(j,j+1)} &=&
    + \ q_{12}^{-1} c_{\uparrow j}^+ (1 - n_{\downarrow j})
       c_{\uparrow {j+1}}^- (1 - n_{\downarrow {j+1}})
    -  q_{12}\ c_{\uparrow j}^- (1 - n_{\downarrow j})
       c_{\uparrow{j+1}}^+ (1 - n_{\downarrow {j+1}}) \nonumber \\
& & +\ q_{23}\ c_{\downarrow j}^+ (1 - n_{\uparrow j})
       c_{\downarrow {j+1}}^- (1 - n_{\uparrow {j+1}})
    -  q_{23}^{-1} c_{\downarrow j}^- (1 - n_{\uparrow j})
       c_{\downarrow{j+1}}^+ (1 - n_{\uparrow {j+1}}) \nonumber \\
& & +\ q_{13}^{-1} \sigma_j^+ \sigma_{j+1}^-
    +  q_{13}\ \sigma_j^- \sigma_{j+1}^+ \nonumber \\
& & +\ q^{-1}(1 - n_{\downarrow j})(1 - n_{\uparrow {j+1}})
    +  q (1 - n_{\uparrow j})(1 - n_{\downarrow {j+1}}) \;.
\label{HamiltJqij}
\ena
The parameters $q_{ij}$ have the following interpretation: $q_{12}$
(resp. $q_{23}$)
corresponds to a left-right anisotropy of the hopping of the fermionic
state $\uparrow$
(resp. $\downarrow$), whereas $q_{13}$ is an
anisotropy for the magnetic interaction.
\\
For open boundary conditions, it is easy to find a transformation that
eliminates two anisotropy parameters. This transformation is
analoguous to the local rescaling that eliminates the anisotropy of
magnetic interaction in the XXZ open chain.
\\
The algebra defined by the relations (\ref{RLeL}), with the $R$-matrix
of (\ref{Rqij}), can be brought in the two parameter form of
Sect. \ref{quantum II}
by rescaling the $\uL^{\pm}$ matrices by the group-like
element $Z=q_{12}^{H_1}q_{13}^{-H_1-H_2}q_{23}^{-H_2}$.

\subsubsection{$\uqs$-invariant Hamiltonians from the distinguished basis}

\indent

We now use the distinguished basis and its natural Hopf structure
(\ref{Deltadist}) (or
equivalently the fermionic basis and the alternative coproduct
(\ref{Deltafermtilde})) to
construct a quantum chain Hamiltonian.

All the Casimir operators again lead to a unique Hamiltonian (up to
normalization)
\be
\cH_{dist} =
\left(
\begin{array}{ccc|ccc|ccc}
0 & 0 & 0 & 0 & 0 & 0 & 0 & 0 & 0 \\
0 & q & 0 & s^{-1} & 0 & 0 & 0 & 0 & 0 \\
0 & 0 & q & 0 & 0 & 0 & -s^{-1} & 0 & 0 \\
\hline
0 & s & 0 & q^{-1} & 0 & 0 & 0 & 0 & 0 \\
0 & 0 & 0 & 0 & q + q^{-1} & 0 & 0 & 0 & 0 \\
0 & 0 & 0 & 0 & 0 & q^{-1} & 0 & s^{-1} & 0 \\
\hline
0 & 0 & -s & 0 & 0 & 0 & q^{-1} & 0 & 0 \\
0 & 0 & 0 & 0 & 0 & s & 0 & q & 0 \\
0 & 0 & 0 & 0 & 0 & 0 & 0 & 0 & 0
\end{array}
\right) \;.
\label{HamilQuantdist}
\ee
In terms of the fermionic creation and
annihilation operators, one has
\bea
\cH_{dist}^{(j,j+1)} &=& + \ s^{-1} c_{\uparrow j}^+ (1 - n_{\downarrow j})
       c_{\uparrow {j+1}}^- (1 - n_{\downarrow {j+1}})
    -  s\ c_{\uparrow j}^- (1 - n_{\downarrow j})
       c_{\uparrow{j+1}}^+ (1 - n_{\downarrow {j+1}}) \nonumber \\
& & +\ s\ c_{\downarrow j}^+ (1 - n_{\uparrow j})
       c_{\downarrow {j+1}}^- (1 - n_{\uparrow {j+1}})
    -  s^{-1} c_{\downarrow j}^- (1 - n_{\uparrow j})
       c_{\downarrow{j+1}}^+ (1 - n_{\uparrow {j+1}}) \nonumber \\
& & -\ s^{-1} \sigma_j^+ \sigma_{j+1}^-
    - s\ \sigma_j^- \sigma_{j+1}^+ \nonumber \\
& & +\ q\ n_{\emptyset j+1}
    +\ q^{-1} n_{\emptyset j}
    +\ q\ n_{\uparrow j} n_{\downarrow {j+1}}
    +\ q^{-1} n_{\downarrow j} n_{\uparrow {j+1}}\;.
\label{HamiltJdist}
\ena

For $s=1$, this Hamiltonian reduces to the one found in \cite{FeiYue},
whatever the Casimir operator we use for its construction.

\subsubsection{Equivalence with one parameter Hamiltonians}

\indent

The Hamiltonian (\ref{HamiltJQuant}) and the slightly more general
Hamiltonian (\ref{HamiltJqij}) are, for open boundary conditions,
equivalent to the one parameter Hamiltonian one gets by setting
$q_{ij}=1$ in (\ref{HamiltJqij}).
The similarity transformation
\be
\cH \longrightarrow \cO^{-1} \cH \cO
\ee
with
\be
\cO = q_{12}^{-\left( \sum\limits_{i<j} n_{\uparrow i} n_{\emptyset j} \right)}
      q_{23}^{-\left( \sum\limits_{i<j} n_{\emptyset i} n_{\downarrow j}
\right)}
      q_{13}^{-\left( \sum\limits_{i<j} n_{\uparrow i} n_{\downarrow j}
\right)}
\ee
indeed
brings the Hamiltonian $\cH(q,q_{12},q_{13},q_{23})$ to the Hamiltonian
$\cH(q,q'_{12}=1,q'_{13}=1,q'_{23}=1)$.
An analogous transformation was also (and earlier) found in
\cite{Dahmen} in the context of Reaction-Diffusion processes.
\\
The same similarity transformation applied to (\ref{HamilQuantdist})
with $q_{12}=-q_{13}=q_{23}=s$ transforms all its non diagonal
terms to $1$, leading to the Perk--Schultz
Hamiltonian $\cH^{P,M}$ with $P=2$, $M=1$ \cite{PerkS}.

Hence, for spin chains with open boundary conditions, the only
relevant parameter is $q$.

\indent

After use of this similarity transformation on the Hamiltonians
(\ref{HamilQuant}) and (\ref{HamilQuantdist}), they only differ by
$q\leftrightarrow q^{-1}$ (reflection of the chain) and the exchange of
two diagonal terms. For two or three sites, they are obviously
equivalent.

We also proved numerically that their spectra are identical up to
seven sites, which strongly suggests that the $L$-site Hamiltonians are
actually equivalent.

\indent

{}From the properties of the $\hR$-matrix (Eqs. (\ref{QYBEh}) and
(\ref{Rhcheqn})), we deduce that the translated Hamiltonians
$U\equiv \cH - q^{\pm 1}$
satisfy the Hecke algebra
\bea
&& U_{i,i+1}^2 \pm \lambda U_{i,i+1} - 1 = 0\;, \nonumber\\
&& U_{i,i+1}U_{i+1,i+2}U_{i,i+1} = U_{i+1,i+2}U_{i,i+1}U_{i+1,i+2}\;,
\nonumber\\
&& [U_{i,i+1},U_{j,j+1}] = 0 \qquad \hbox{for} \qquad |i-j|\ge 2 \;.
\ena

\subsubsection{Uniqueness of the invariant Hamiltonians}

\indent

Finally, it is interesting to note that one can compute all the possible
$\uqs$-invariant $3^L$-state
Hamiltonians with nearest neighbour interaction directly, by
imposing that
\be
\left[\cH_{12} , X_{12} \right] = 0
\ee
for all the generators $X$ of $\uqs$. The unique solution
up to the identity and a normalization is given by
Eq. (\ref{HamilQuant}) in the fermionic case, and by
Eq. (\ref{HamilQuantdist}) in the distinguished case.
This means that the commutant of $\Delta(\uqs)$ in the tensor square
of the fundamental representation is generated by $I$ and
the matrix $\hR$ (the matrix $\hR$ being either the fermionic or the
distinguished one).

The physical consequence is that there are two  deformations of
the $t-J$ model with $\uqs$-invariance, actually corresponding to the
two different Hopf structures (\ref{Deltaferm}) and
(\ref{Deltafermtilde}) on $\uqs$.

\sect{Summary \label{Conclusion}}

\indent

In the case of classical $sl(1|2)$, we found that the centre is
generated
by an infinite number of Casimir operators $\cC_k$ ($k$ positive
integer)
of degree $k$ in the Cartan generators. However, the $\cC_k$ are not
completely independent since they obey quadratic relations. It follows
that any two Casimir operators are sufficient to label an irreducible
finite dimensional representation of $sl(1|2)$.
Using the construction of Sect. \ref{Hamil},
one can show that all the Casimir operators
$\cC_k$ lead to the same two-site Hamiltonian,
up to the identity matrix and a normalization
factor. Moreover,  the same results hold if one uses either
the fermionic basis or the distinguished basis of $sl(1|2)$, i.e. the
two simple root bases of $sl(1|2)$.

In the deformed case, the situation is more complicated.

In the two-parametric deformed case, the centre of $\uqs$ is
generated by an infinite number of Casimir operators $\cC_k$ ($k \in \ZZ$),
which obey the same quadratic relations as in the classical case.
In the fermionic basis, the Casimir operators $\cC_k$ lead to two-site
Hamiltonians which are all proportional to the same $\cH$.

Such Hamiltonians describe deformed $t-J$ models at the
supersymmetric point, the deformation parameter $s$ being interpreted
as a left-right anisotropy for the hopping terms and the magnetic
interaction and $q$ as an anisotropy in $sl(1|2)$ analoguous to that
of the usual XXZ chain.
The deformation parameter $s$ (and more generally the
parameters $q_{ij}$ of Eq. (\ref{Hamilqij})) can be removed, for an
open spin chain,  by a
similarity transformation.

In the distinguished basis, all the Casimir operators lead to a unique
Hamiltonian, which is the one found in \cite{FeiYue}.
It is also a deformation of the $t-J$ model at the
supersymmetric point, and for $s=1$, corresponds to the
Perk--Schultz $(1,2)$ Hamiltonian.
We checked numerically that both Hamiltonians
(\ref{HamilQuant}) and (\ref{HamilQuantdist}) are actually equivalent.

Lets us finally note that using the conjugate of the fundamental
representation amounts to a reversal of the spin chain.

\newpage

\sect{Appendix \label{Appendix}}

\indent

We give here the explicit relations among the $L$ matrices. In these
relations, the value $\zeta=1$ corresponds to the $Z_2$-graded case
($\uL$ matrices of Section \ref{quantum II}),
whereas the value $\zeta=-1$ corresponds to the bosonised case
($L$ matrices of Section \ref{Boson}). For simplicity, we
use a single notation $L$ in this Appendix.
\begin{eqnarray*}
&&
\begin{array}{ll}
  L^{\pm}_{11} L^+_{12} = \zeta q^{\mp 1} s^{-1} L^+_{12} L^{\pm}_{11} \;,
& L^{\pm}_{11} L^-_{21} = \zeta q^{\pm 1} s L^-_{21} L^{\pm}_{11} \;, \\
  L^{\pm}_{22} L^+_{12} = q^{\mp 1} s^{-1} L^+_{12} L^{\pm}_{22} \;,
& L^{\pm}_{22} L^-_{21} = q^{\pm 1} s L^-_{21} L^{\pm}_{22} \;, \\
  L^{\pm}_{33} L^+_{12} = \zeta L^+_{12} L^{\pm}_{33} \;,
& L^{\pm}_{33} L^-_{21} = \zeta L^-_{21} L^{\pm}_{33} \;, \\
& \\
  L^{\pm}_{11} L^+_{23} = \zeta L^+_{23} L^{\pm}_{11} \;,
& L^{\pm}_{11} L^-_{32} = \zeta L^-_{32} L^{\pm}_{11} \;, \\
  L^{\pm}_{22} L^+_{23} = q^{\pm 1} s^{-1} L^+_{23} L^{\pm}_{22} \;,
& L^{\pm}_{22} L^-_{32} = q^{\mp 1} s L^-_{32} L^{\pm}_{22} \;, \\
  L^{\pm}_{33} L^+_{23} = \zeta q^{\pm 1} s^{-1} L^+_{23} L^{\pm}_{33} \;,
& L^{\pm}_{33} L^-_{32} = \zeta q^{\mp 1} s L^-_{32} L^{\pm}_{33} \;, \\
& \\
  L^{\pm}_{11} L^+_{13} = q^{\mp 1} s^{-1} L^+_{13} L^{\pm}_{11} \;,
& L^{\pm}_{11} L^-_{31} = q^{\pm 1} s L^-_{31} L^{\pm}_{11} \;, \\
  L^{\pm}_{22} L^+_{13} = s^{-2} L^+_{13} L^{\pm}_{22} \;,
& L^{\pm}_{22} L^-_{31} = s^2 L^-_{31} L^{\pm}_{22} \;, \\
  L^{\pm}_{33} L^+_{13} = q^{\pm 1} s^{-1} L^+_{13} L^{\pm}_{33} \;,
& L^{\pm}_{33} L^-_{31} = q^{\mp 1} s L^-_{31} L^{\pm}_{33} \;, \\
& \\
  L^+_{12} L^+_{12} = 0 \;,
& L^-_{21} L^-_{21} = 0 \;, \\
  L^+_{23} L^+_{23} = 0 \;,
& L^-_{32} L^-_{32} = 0 \;, \\
& \\
  L^+_{12} L^+_{23} + \zeta L^+_{23} L^+_{12} = \lambda s^{-1} L^+_{13}
L^+_{22} \;,
& L^-_{32} L^-_{21} + \zeta L^-_{21} L^-_{32} = - \lambda s L^-_{31} L^-_{22}
\;,  \\
  L^+_{13} L^+_{12} = \zeta q s L^+_{12} L^+_{13} \;,
& L^-_{31} L^-_{32} = \zeta q^{-1} s^{-1} L^-_{32} L^-_{31} \;,  \\
  L^+_{23} L^+_{13} = \zeta q s^{-1} L^+_{13} L^+_{23}  \;,
& L^-_{31} L^-_{21} = \zeta q s^{-1} L^-_{21} L^-_{31} \;, \\
& \\
  L^+_{12} L^-_{32} = -q^{-1} s L^-_{32} L^+_{12} \;,
& L^+_{23} L^-_{21} = -q s L^-_{21} L^+_{23} \;, \\
& \\
  s^{-1} L^+_{12} L^-_{31} - \zeta s L^-_{31} L^+_{12}
     = \lambda L^-_{32} L^+_{11} \;, \\
  s^{-1} L^+_{13} L^-_{21} - \zeta s L^-_{21} L^+_{13}
     = - \zeta \lambda L^+_{23} L^-_{11} \;, \\
  s^{-1} L^+_{23} L^-_{31} - \zeta s L^-_{31} L^+_{23}
     = - \lambda L^-_{21} L^+_{33} \;, \\
  s^{-1} L^+_{13} L^-_{32} - \zeta s L^-_{32} L^+_{13}
     = \zeta \lambda L^+_{12} L^-_{33} \;, \\
& \\
  s^{-1} L^+_{12} L^-_{21} + \zeta s L^-_{21} L^+_{12}
     = - \lambda (L^+_{11} L^-_{22} - L^+_{22} L^-_{11}) \;, \\
  s^{-1} L^+_{23} L^-_{32} + \zeta s L^-_{32} L^+_{23}
     = \zeta \lambda (L^+_{22} L^-_{33} - L^+_{33} L^-_{22}) \;, \\
  s^{-1} L^+_{13} L^-_{31} - s L^-_{31} L^+_{13}
     = \lambda (L^+_{11} L^-_{33} - L^+_{33} L^-_{11}) \;.
\end{array}
\end{eqnarray*}

\indent

{\bf Acknowledgements:}

\medskip

We would like to acknowledge M. Bauer, S. Dahmen, E. Ragoucy,
V. Rittenberg and P. Sorba for fruitful discussions and for
communicating some useful references and papers.


\end{document}